\documentclass[]{article}

\usepackage{amsmath, amssymb, amsfonts, amsthm, latexsym}
\usepackage{graphicx}

\theoremstyle{plain}
\newtheorem{theorem}{Theorem}[section]

\newtheorem{lemma}[theorem]{Lemma}
\newtheorem{prop}[theorem]{Proposition}
\newtheorem{cor}[theorem]{Corollary}

\theoremstyle{definition}
\newtheorem{defi}[theorem]{Definition}
\newtheorem{note}[theorem]{Notation}
\newtheorem{ex}[theorem]{Example}

\newtheorem{rmk}[theorem]{Remark}
\newtheorem{remark}[theorem]{Remark}

\newcommand{\Gbb}{{\mathbb G}}

\newcommand{\codim}{{\rm codim\ }}

\newcommand{\Gcal}{{\mathcal G}}

\newcommand{\Ical}{{\mathcal I}}
\newcommand{\Jcal}{{\mathcal J}}
\newcommand{\Bcal}{{\mathcal B}}

\newcommand{\Kcal}{{\mathcal K}}
\newcommand{\Ccal}{{\mathcal C}}
\newcommand{\Rcal}{{\mathcal R}}

\newcommand{\la}{\langle}
\newcommand{\ra}{\rangle}

\newcommand{\btu}{\bigtriangleup}

\newcommand{\f}[6]{x^{#1} x^{#2} x^{#3} - x^{#4} x^{#5} x^{#6}}

\newcommand{\abs}[1]{\mid\!\!#1\!\!\mid}
\newcommand{\ed}[2]{\langle #1, \ #2 \rangle}
\newcommand{\fr}[1]{\mathfrak #1}

\newcommand{\ov}[1]{\overline{#1}}
\newcommand{\wt}[1]{\widetilde {#1}}

\author{HA Minh Lam\\ MORALES Marcel}
\begin{document}
\title{Binomial extensions of Simplicial  ideals and reduction number}
\date{}
\maketitle
\section*{Abtracts}
In this article, we define a class of binomial ideals associated to a simplicial complex.  This class of ideals appears in the presentation of fiber cones of codimension 2 lattice ideals \cite{hm}, and in the work of Barile and Morales \cite{bm2}, \cite{bm3}, \cite{bm4}.
We  compute the reduction number of Binomial extensions of Simplicial  ideals. This extends all the previous results in this area.
\section*{Introduction}

According to the classification resulting from the successive contributions
by Del Pezzo, Bertini, and Xamb\'o (see \cite{eg} for literature), the equidimensional algebraic subsets $X\subset \mathbb P^n $ of minimal degree which are connected in codimension one are of three types: quadric hypersurfaces, the cone over the Veronese surface
in $\mathbb P^5$, and unions $X = \cup^n_{i=0} X_i$ of scrolls embedded in linear subspaces such that for all $i = 1, \ldots, n - 1$, we have:
$$(X_1 \cup \ldots \cup X_i) \cap X_{i+1} = {\rm Span}(X_1 \cup \ldots \cup X_i) \cap  {\rm Span}(X_{i+1}). \qquad  (*)$$ 
Homologically, varieties of minimal degree were characterized by Eisenbud --
Goto \cite{eg} (see Theorem \ref{Eisenbud-Goto}).

Under the algebraic point of view, the condition $(*)$ was considered at first in \cite{bm2}, and later in \cite{bm3}, the authors  give a complete constructive characterization of the ideals defining varieties of unions of scrolls satisfying the above condition $(*)$.

Later Eisenbud--Green--Hulek--Popescu \cite{eghp}, define a {\it linearly joined} sequence of varieties, as an union of varieties satisfying the condition $(*)$.  They prove that an algebraic set $X \subset \mathbb P^r$ is 2--regular if, and only if, $X = X_1 \cup \ldots X_n,$ with $X_1, \ldots ,X_n$ is a  sequence of varieties of minimal degree.

Recall that the homogeneous coordinate ring of a scroll is of the type $A = S/I_2(M)$, where $S$ is the polynomial ring
$S = k[{T_{i,j} \mid 1 \le i \le l\ \rm{and}\ 1 \le j \le s_i + 1}]$, and $I_2(M)$ is the ideal generated by $2 \times 2$ minors of the matrix $M$ = $( M_1\mid M_2\mid \ldots\mid M_l)$, with each $M_u$ is the generic catalecticant
matrix
$$M_u = \left(\begin{array}{ccccc}T_{u,1}& T_{u,2}& \ldots &T_{u,s_u-1} &T_{u,s_u}\\T_{u,2}&T_{u,3}& \ldots  &T_{u,s_u}&T_{u,s_u+1}\end{array}\right).$$
We call an ideal of type $I_2(M)$ a scroll ideal.

In this article, we define a class of binomial ideals associated to a simplicial complex.  This class of ideals appears in the presentation of fiber cones of codimension 2 lattice ideals \cite{hm}, and in the work of Barile and Morales \cite{bm2}, \cite{bm3}, \cite{bm4}.
Let $\btu$ be a simplicial complex over a set of vertices $V_\bigtriangleup = \{ x_1, x_2, \ldots , x_n \}.$ We will call proper facet a facet $F_l$ with a star of some edges belonging only to $F_l$ (called also  {\it proper edges}). To each {\it proper facet} of $\btu$,  we associate a set of points $Y^{(l)}$ (which can be empty), and a scroll ideal $I_l$ of variables in $Y^{(l)}$ and in vertex set of these proper edges. The new simplicial complex obtained from $\btu$ and the sets $Y^{(l)}$ is called an {\it extension complex}, and denoted by $\ov \btu$.

The {\it binomial extension of a simplicial  ideal} $\Bcal_{\overline{\bigtriangleup}}$ associated to $\ov \btu$ is defined to be the one generated by all $I_l$ and the Stanley--Reisner ideal of $\ov \btu$.

The aim of this article is to prove that binomial extension of simplicial  ideal is a good generalization of Stanley--Reisner theory to the case of binomial ideals. 

In the first section, we will define the class of binomial extension of simplicial  ideals, we will give  the prime decomposition. From this, we deduce that our class of ideals, in fact, defines an union of scrolls along linear spaces.

In Section 2, we study the reduction number of binomial extension of  simplicial  ideals. Our aim is to extend the results of Barile and Morales : to describe explicitly the reduction ideals through the complexes. In \cite{bm1}, they described a class of square--free monomial ideals whose reduction number is 1 by coloring the graph of a simplicial complex.
In \cite{bm2}, Barile and Morales considered a class of binomial ideals, which indeed are particular cases of binomial extension of  simplicial  ideals $\Bcal_{\ov \btu}$ where $\Gcal_{\btu}$ is a generalized $d$--tree and each vertex belongs to at most two extension facets. They proved that this class of ideals is of reduction number 1, and an explicit expression of the reduction is given. 

In the case of binomial extension of  simplicial  ideals, we have the following theorem:
\begin{theorem} If the graph associated to $\btu$ admits a good $(d+1)$--coloration, and in addition, for each proper facet $F$ the origin of the star of proper edges belongs only to $F$, then the ring $K[{\rm {\bf x}}, {\rm {\bf y}}] \slash \Bcal_{\ov \btu}$ has reduction number 1. 

In this case, the reduced graph associated to $\ov \btu$ admits also a good $(d+1)$--coloration, and
$$(g_1,\ g_2, \ldots, \ g_{d+1}){\fr m}_{\ov \btu} + \Bcal_{\ov \btu} = \fr m^2_{\ov \btu},$$
where ${\fr m}_{\ov \btu} = \left({\rm {\bf x}}, {\rm {\bf y}} \right)$ is the irrelevant ideal of $K[{\rm {\bf x}}, {\rm {\bf y}} ],$ and $g_i$ is the sum of all variables with color $i$.
\end{theorem}

\begin{theorem} Let $\Gcal_{\btu}$ be a generalized $d$--tree. Then we can find a good $(d+1)$--coloration for the reduced graph associated to $\ov \btu$, such that
$$(g_1,\ g_2, \ldots, \ g_{d+1}){\fr m}_{\ov \btu} + \Bcal_{\ov \btu} = \fr m^2_{\ov \btu},$$
where ${\fr m}_{\ov \btu} = \left({\rm {\bf x}}, {\rm {\bf y}} \right)$ is the irrelevant maximal ideal of the polynomial ring $K[{\rm {\bf x}}, {\rm {\bf y}} ],$ and $g_i$ is the sum of all variables with color $i$.
\end{theorem}

% and their result offers a way to generalize the hypothesis of the Del
%Pezzo-Bertini Theorem to all projective schemes. To state their result recall that
%$X \subset \Pbb^r$ is said to be $d$-regular, in the sense of Castelnuovo-
%Mumford, if the ideal sheaf $\Ical_X$ satisfies $H_i(\Ical_X(d - i)) = 0$ for all $i > 0$, or equivalently,
%if the $j^{\rm th}$ syzygies of the homogeneous ideal $\Ical_X$ are generated in degrees $\le d + j$ for all $j \ge 0$ (see Eisenbud-Goto \cite{eg} for a proof of the equivalence). 

\section{Simplicial ideals and  binomial extension of  simplicial  ideals }

binomial extension of  simplicial  ideals is an extension of Stanley-Reisner monomial ideals. It associates an ideal to a simplicial complex and a family of ideals indexed by a set of its facets. In this article, we will consider some particular cases, which define in fact an union of scrolls, and study some properties of these binomial extension of  simplicial  ideals.

First of all, let us recall some definitions.

A \emph{simplicial complex} $\bigtriangleup$ over a vertex set $V_\bigtriangleup = \{x_1, x_2, \ldots , x_n\}$ is a collection of subsets of $V_\bigtriangleup$ with the property that:
\begin{itemize} \item For all $i,$ the set $\{x_i\}$ is in $\bigtriangleup,$
\item If $F \in \bigtriangleup$ and $G \subset F,$ then $G \in \bigtriangleup.$
\end{itemize}

An element of a simplicial complex $\bigtriangleup$ is called a \emph{face} of $\bigtriangleup.$ The dimension of a face $F$ of $\bigtriangleup,$ denoted by $\dim F,$ is defined to be $\abs F - 1,$ where $\abs F$ denotes the number of vertices in $F.$ The dimension of $\bigtriangleup,$ denoted by $\dim \bigtriangleup,$ is defined to be the maximal dimension of a face in $\bigtriangleup.$ The maximal faces of $\bigtriangleup$ under inclusion are called \emph{facets} of $\bigtriangleup.$

Let us remark that by taking all faces of dimension $0$ and $1$ of $\bigtriangleup,$ i.e. all vertices and edges, we associate to $\bigtriangleup$ a simple graph $\Gcal_{\btu}.$ An arbitrary facet $F$ of $\bigtriangleup$ becomes a completed subgraph of $\Gcal_\btu$, so called a $\abs F$-clique of $\Gcal_\btu.$
\begin{note} We denote by $(E)$ (resp. $\Kcal [E]$) the ideal generated by $E$ (resp. the polynomial ring with the variables in $E$). And for a facet $F$, let $F^\circ$ denote the set of the points of $F$ which belong to no other facets of $\btu.$  To a set of vertices $V_\bigtriangleup = \{ x_1, x_2, \ldots , x_n \},$ one associate a ring of polynomials $R = \Kcal[x_1, x_2, \ldots , x_n]$ (here, by abuse of notation, we use the $x_i$'s to denote both the vertices in $V_\bigtriangleup$ and the variables in the polynomial ring). It is known that to a simplicial complex $\bigtriangleup$ on this vertex set, there is an associated ideal, called the Stanley -- Reisner ideal, defined as follows
$$\Ical_\bigtriangleup = \left(x_{i_1}x_{i_2} \cdots x_{i_r} \mid i_1 < i_2 < \cdots < i_r, \{x_{i_1}, x_{i_2}, \ldots, x_{i_r}\} \notin \bigtriangleup \right). $$
This ideal is generated by monomials.
\end{note} 
Now, we introduce the definition of a binomial extension of a simplicial complex.   
\begin{defi}A facet $F_l$ of dimension $d_l$ of $\bigtriangleup$ is said to be \emph{proper} if $F$ contains $k_l$ edges $\ed {x^{(l)}_0} {x^{(l)}_{i_1}}, \ed {x^{(l)}_0} {x^{(l)}_{i_2}}, \ldots, \ed {x^{(l)}_0} {x^{(l)}_{i_{k_l}}}$ which belong uniquely to $F$. If it is the case, these edges are called the \emph{proper edges} of $F$.

To each proper edge $\ed {x^{(l)}_0} {x^{(l)}_{i_j}}$ of a  facet $F_l $ of $\bigtriangleup,$ we associate a set of points $Y^{(l)}_j$ (which can be empty). The simplex which is the product of $Y^{(l)}:= \cup Y^{(l)}_j$ and $F_l$ is called {\it the extension of } $F_l$ by $Y^{(l)}$.  By abuse of notation, we use $F_l$ to denote this  new facet (without any confusion). Let us remark that this extension can be trivial for some facets of $\btu$. 
Let $\ov \btu$ be the complex which facets are all the extension facets of $\btu$. We call this complex an {\it extension complex } of $\btu$.
\end{defi}
Let $\overline{\bigtriangleup}$ be an extension complex constructed by $\bigtriangleup$ and a set of points. 
We associate to $\overline{\bigtriangleup}$ a polynomial ring $\mathcal R := \Kcal[V_{\overline{\bigtriangleup}}]=\Kcal[{\rm {\bf x}}, {\rm {\bf y}}]$. We will denote by $x^{(l)}_{i_j}$ the vertices in $F_l\cap \bigtriangleup$   and by $y^{(l)}_{ij_m}$ the vertices in  $F_l\setminus \bigtriangleup$.

\begin{center}
\includegraphics[height=7.3cm,width=8.4cm]{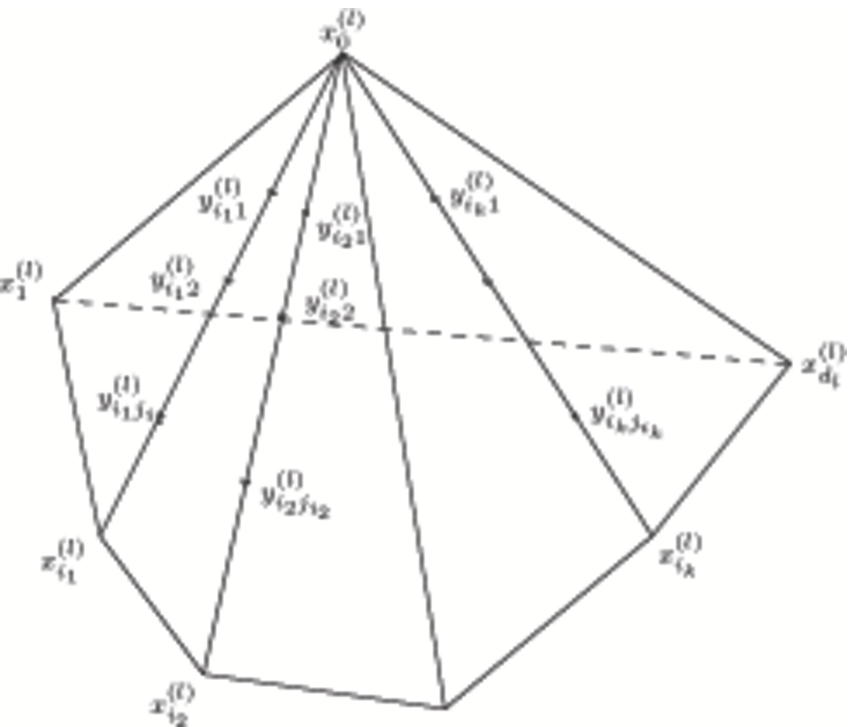}
\end{center}

\begin{defi}To each non trivial extension facet $F_l$ of $\ov \btu$ we associate the prime ideals $\Ical_l, \Jcal_l$, where
\begin{description} \item $\bullet$ $\Ical_l = 0$ if $Y^{(l)} = \emptyset$, 
\item $\bullet$ if $Y^{(l)} := \{y^{(l)}_{i_jm}\} \neq \emptyset$, the $\Ical_l$ is the ideal generated by the $2 \times 2$ minors of the matrix:
$$ M_l:= \left( \begin{array}{cccc} x^{(l)}_0 & y^{(l)}_{i_11} & \ldots & y^{(l)}_{i_1j_{i_1}}\\ y^{(l)}_{i_11} 
& y^{(l)}_{i_12} & \ldots & x^{(l)}_{i_1} \end{array} \left|
\begin{array}{ccc} y^{(l)}_{i_21} & \ldots & y^{(l)}_{i_2j_{i_2}} \\ y^{(l)}_{i_22} & \ldots & x^{(l)}_{i_2}
\end{array} \left|
\begin{array}{c} \ldots \\ \ldots \end{array} \left|
\begin{array}{ccc}  y^{(l)}_{i_{k_l}1} & \ldots & y^{(l)}_{i_{k_l}j_{i_{k_l}}}\\ y^{(l)}_{i_{k_l}2} &\ldots &
 x^{(l)}_{i_{k_l}}   
    \end{array}   \right. \right. \right.      \right).$$
\item  $\Jcal_l= \left(\Ical_l, (V_{\ov \btu}\setminus F_l) \right)\subset \Kcal[V_{\overline{\bigtriangleup}}],$ where $V_{\ov \btu}\setminus F_l$ denote the set of vertices of $\ov \btu$ which are not in $ F_l.$ 
\end{description}
The {\it binomial extension of  simplicial  ideal} $\Bcal_{\overline{\bigtriangleup}}\subset \Kcal[V_{\overline{\bigtriangleup}}]$ is defined by 
$$\Bcal_{\ov \btu} = \left( \sum_{F_l \ \rm{ facet \ of }\ \ov \btu} \Ical_l, \Ical_{\ov \btu} \right),$$
where $\Ical_{\ov \btu}$ is the Stanley--Reisner ideal associated to the simplicial complex $\overline{\bigtriangleup}.$

The couple $\btu(\Bcal_{\ov \btu}) := ({\ov \btu}, \Bcal_{\ov \btu})$ is called a {\it binomial extension} of $\btu$. 
\end{defi}

It is well--known that the Stanley--Reisner ideal of a simplicial complex admits a decomposition into prime ideals corresponding to the facets of the complex (each ideal is generated by the variables which are not in the correspondent facet). We will prove the same property for the ideal $\Bcal_{\ov \btu}$.

%\begin{itemize}
%\item $\Jcal_l =\left(V_{\ov \btu}\setminus F_l \right),$ if $F_l$ is also a facet of $\btu$ (in this case $\Ical_l=0$).
%\item $\Jcal_l = \left(\Ical_l, (V_{\ov \btu}\setminus F_l) \right),$ otherwise
%\end{itemize}
%where $V_{\ov \btu}\setminus F_l$ denote the set of vertices of $\ov \btu$ which are not in $ F_l.$
 
We have the Primary decomposition of $\Bcal_{\overline{\bigtriangleup}}$ in the ring $\Kcal[V_{\overline{\bigtriangleup}}]$:
\begin{prop}\label{ibin} $$\displaystyle \Bcal_{\overline{\bigtriangleup}} = 
\bigcap_{F_l \mbox{ facet of } {\ov \bigtriangleup}} \Jcal_l.$$
\end{prop}

Before proving the proposition, let us remark that \begin{itemize}
\item $x^{(l)}_0 x^{(l)}_{i_{s}}\in (V_{\overline{\bigtriangleup}}\setminus F_i)$ for all proper edge $\langle{x^{(l)}_0},{x^{(l)}_{i_{s}}}\rangle$ of $F_l$ $(s = \overline{1,k_l})$, and all facet $F_i$ of $\overline{\bigtriangleup}$, $F_i \neq F_l$, since either $x^{(l)}_0 \in V_{\overline{\bigtriangleup}}\setminus F_i$ or $x^{(l)}_{i_{s}}\in V_{\overline{\bigtriangleup}}\setminus F_i$.
\item $(F_l^\circ)\subset (V_{\overline{\bigtriangleup}}\setminus F_i)$, for all facet $F_i$ of $\overline{\bigtriangleup}$, $F_i \neq F_l$.
\end{itemize} 
It implies that for all facet $F_i$ of $\overline{\bigtriangleup}$, $F_i \neq F_l$ we have 
$$\Ical_l \subset (V_{\overline{\bigtriangleup}}\setminus F_i) \eqno (*)$$ 

{\sc Proof:}  First of all, recall a well-known fact that if $I, J$ are disjoint sets of variables and 
 $\Ical  = \left( I  \right), \Jcal  = \left( J \right)$ then
 $\Ical\cap  \Jcal = \left(p q\left|\  p \in I,q\in J \right. \right).$

We will prove the proposition by induction on the number $m$ of facets of $\ov \btu.$ The case $m=1$ is trivial. If $m >1$, then we denote by $F_m $ the $m^{\rm{th}}$ facet. Denote by $\ov {\btu^\prime}$ the complex constructed by $m-1$ facets of $\ov \btu$, and denote by $\Bcal_{\ov {\btu^\prime}} $ the binomial ideal associated to $\ov {\btu^\prime}$, and $\Jcal^\prime_i$ the prime ideal associated to the $i^{\rm{th}}$ facet of $\ov {\btu^\prime}.$
 Remark that $\Jcal_i = \left(\Jcal^\prime_i,\ x \left|\ x \in F_m^\circ \right) \right.$ for all $(i = 1, 2,\dots, m-1).$ By induction, we have $\Bcal_{\ov {\btu^\prime}}= \bigcap \Jcal^\prime_i.$ 
We have that $\Jcal_m = \left( \Ical_m, V_{\ov \btu}\setminus F_m \right),$ and
$$\Bcal_{\ov {\btu^\prime}}    = \left(\sum_{1\leq k\leq m-1} \Ical_k, \Ical_{\ov {\btu^\prime}} \right), \mbox{ and  }\Bcal_{\ov \btu} = \left(\sum_{1\leq k\leq m} \Ical_k + \Ical_m, \Ical_{\ov \btu} \right).$$
For the Stanley--Reisner ideal $\Ical_{\ov \btu}$, it is known that:
$$ \left(\Ical_{\ov {\btu^\prime}} \ , \ F_m^\circ  \right) \bigcap \left( V_{\ov \btu}\setminus F_m  \right) 
= \Ical_{\ov \btu}. \eqno (\alpha)$$
In addition, since $\Ical_k \subset \left( V_{\ov \btu}\setminus F_m\right)$ for all $k \neq m$, we have: 
 $$ \sum_{1\leq k\leq m-1} \Ical_k \subset \left( V_{\ov \btu}\setminus F_m\right).
 \eqno (\beta)$$
Moreover, since $\Ical_m \subset \left( V_{\ov \btu}\setminus F_{m^\prime}  \right)$
for all facet $F_{m^\prime} \neq F_m$ in $\ov \btu$, we have:
$$\Ical_m \subset \bigcap_{{m^\prime} \neq m} \left( V_{\ov \btu}\setminus F_{m^\prime} \right) =
 \left(\Ical_{\ov {\btu^\prime}},  F_m^\circ  \right). \eqno (\gamma)$$
It implies that:
$$\Bcal_{\ov \btu} \subset \left(\Bcal_{\ov {\btu^\prime}},\ F_m^\circ  \right) \bigcap 
\left( \Ical_m, V_{\ov \btu}\setminus F_{m} \right).$$
Now, we will prove the other inclusion. If $r \in \left(\Bcal_{\ov {\btu^\prime}},\ F_m^\circ  \right) \bigcap 
\left( \Ical_m, V_{\ov \btu}\setminus F_{m} \right),$ then $r=u+v = p + q$, where
$u\in \sum_{1\leq k\leq m-1} \Ical_k,$   $v\in \left(\Ical_{\ov {\btu^\prime}},\ F_m^\circ  \right)$, 
 $p \in \Ical_m,$ and $q\in \left(   V_{\ov \btu}\setminus F_{m}  \right)$. 
Due to $(\beta)$ and $(\gamma)$, we have
$$v - p = q - u \in \left(V_{\ov \btu}\setminus F_{m}  \right) \cap \left(\Ical_{\ov {\btu^\prime}},
 F_m^\circ  \right) \stackrel{(\alpha)}{=} \Ical_{\ov \btu}.$$ 
Hence, 
$$r = u + (v-p) + p \in \sum_{1\leq k\leq m-1} \Ical_k + \Ical_{\ov \btu} + \Ical_m = \Bcal_{\ov \btu}.$$
From this it follows that:
$$\Bcal_{\ov \btu} = \left(\Bcal_{\ov {\btu^\prime}},\ F_m^\circ  \right) \bigcap 
\left( \Ical_m, V_{\ov \btu}\setminus F_{m} \right).$$
From the induction hypotheses, we deduce that
$$\Bcal_{\ov \btu} = \left(\bigcap^{m-1}_{i=1} \Jcal^\prime_i, \ F_m^\circ  \right) \bigcap \Jcal_m = 
\bigcap^{m-1}_{i=1} (\Jcal^\prime_i, \ F_m^\circ ) \bigcap \Jcal_m = \bigcap^{m-1}_{i=1} \Jcal_i \bigcap \Jcal_m.$$
The proposition is proved. \hfill $\Box$\\

\begin{rmk} For all facet $F_l$ of $\ov \btu$, the ideal $\Jcal_l$ is prime and the ring $\Kcal[V_{\overline{\bigtriangleup}}]/\Jcal_l$ is of dimension $1 + d_l,$ where $d_l$ is the dimension of $F_l$.
\end{rmk}

We deduce from that a corollary on the dimension of $\Bcal_{\ov \btu}$ as follows:

\begin{cor} $\dim(\Kcal[V_{\overline{\bigtriangleup}}]/\Bcal_{\overline{\bigtriangleup}}) = 1 + \dim(\ov \btu).$
\end{cor}

%Let us remark that the precedent proposition is also valid for a simplicial ideal since the key point $(*)$ is still verified for the simplicial ideals. Recall that 

%\begin{defi} To each facet $F_l,$ we associate a prime ideal $\Ical_l\subset (F_l)$ in the ring $\Kcal [V_{\overline{\bigtriangleup}} ]$ such that $\Ical_l \subset (V_{\overline{\bigtriangleup}}\setminus F_i)$ for all facet $F_i \neq F_l$ of $\overline{\bigtriangleup}$.
 
%The {\it simplicial ideal} $\Pcal_{\ov \btu}$ associated to  the extension complex $\ov \btu$ is defined by
%$$\Pcal_{\ov \btu} = \left( \sum_{F_l \ \rm{ facet \ of }\ \ov \btu} \Ical_l, \Ical_{\ov \btu} \right),$$
%where $\Ical_{\ov \btu}$ is the Stanley--Reisner ideal associated to the simplicial complex $\overline{\bigtriangleup}.$
%\end{defi}

\section{ Reduction number one}

First, we recall a theorem of Eisenbud--Goto \cite{eg}:
\begin{theorem}\label{Eisenbud-Goto}
Let $R$ be a reduced graded ring, defining an algebraic projective variety. Then we have:
\begin{description}
\item{$\quad \ (1)$} $R$ is Cohen-Macaulay and $e(R)=1+\codim R$, where $e(R)$ is the multiplicity of $R$;
\item{$\Rightarrow$ (2)} $R$ admits a $2-$linear resolution;
\item{$\Rightarrow$ (3)} $r(R)=1$;
\item{$\Rightarrow$ (4)} $e(R)\leq 1+\codim R$
\end{description}
Moreover, if $R$ is Cohen--Macaulay, then the above implications are equivalences.
\end{theorem}

\begin{defi} A {\it generalized $d-$tree} on a set of vertices $V$ is a graph defined recursively by the following properties:
\begin{description}
\item{(a)} A complete graph on $d+1$ elements of $V$ is a generalized $d-$tree.
\item{(b)} Let $G$ be a graph on the set $V$. Assume that there exists a vertex $v\in V$ such that: 
\begin{enumerate}
\item The restriction $G'$ of $G$ on $V'=V\setminus\{v\}$ is a generalized $d-$tree,
\item There is a subset $V''\subset V'$ with $1\leq j\leq d$ vertexes such that the restriction of $G$ on $V''$ is a complete graph, and 
\item $G$ is the graph generated by $G'$ and the complete graph on $V''\cup\{v\}$.
\end{enumerate}
 \end{description}
The vertex $v$ as above is called a {\it extremal}. 
 
If $j=d$, then we say that $G$ is a {\it $d-$tree}.
\end{defi}
\begin{rmk}\label{arbre-morales} Let $\Delta (G)$  the ``clique complex" of $G$, i.e. the simplicial complex whose vertices are the ones of $G$ and the facets are the simplexes with support on the complete subgraphs of $G$.  M. Morales associate to $\Delta (G)$ a graph $H(G)$ whose vertices are the facets of $\Delta (G)$, and an edge of $H(G)$ links two vertices such that the intersection of their associated facets is non--empty. He proved that $G$ is a generalized $d-$tree if and only if $H(G)$ is a tree.
\end{rmk}
The following theorems are proved by Fr\"oberg \cite{fr}:
\begin{theorem} The Stanley--Reisner ring of a simplicial complex $\Delta $ is a Cohen--Macaulay ring of minimal degree if and only if 
\begin{enumerate}
\item The graph $ \Gcal_{ \btu}$ is a $d-$tree, and
\item $\Delta $ is a clique complex of $ \Gcal_{ \btu}$, i.e. $ \Delta =\Delta (\Gcal_{ \btu}).$
\end{enumerate}
\end{theorem}
\begin{theorem} The Stanley--Reisner ring of a simplicial complex $\Delta $ admits a $2-$linear resolution if and only if
\begin{enumerate}
\item The graph  $ \Gcal_{ \btu}$  is a generalized $d-$tree, and
\item $\Delta =\Delta (\Gcal_{ \btu}).$
\end{enumerate}
\end{theorem}

\begin{defi} Let $\Rcal$ be the polynomial ring of $n$ variables on the field $\Kcal$. Let $\Ical\subset \Rcal$ be a homogeneous graded ideal under the standard graduation and $d= \dim \Rcal/\Ical.$ A set of linear forms $\{g_1,\ g_2, \ldots, \ g_{d}\}$ is a \emph{reduction} of $\Rcal/\Ical$, if
 $$(g_1,\ g_2, \ldots, \ g_{d}){\fr m}^\rho   = \fr m^{\rho +1} (\mod \Ical) $$
where ${\fr m}$ is the irrelevant maximal ideal of $\Rcal.$

The smallest number $\rho $ for all the possible reductions is called \emph{the reduction number} of $\Rcal \slash \Ical$.
\end{defi}
In \cite{bm1}, Barile and Morales described a class of square--free monomial ideals whose reduction number is 1. First of all, we recall some definitions.
\begin{defi} A \emph{$(d+1)$-coloration} of a graph $\Gcal$ is a
  partition of the vertex set $V_\Gcal$ into $d+1$ subsets, which are called ``class of colors", such that two neighbors in $\Gcal$ belong to different classes of colors. For each vertex $x \in \Gcal,$ we denote by $\Ccal (x)$ the class containing $x.$

A $(d+1)$-coloration of $\Gcal$ is \emph{good} if every cycle of $\Gcal$ is colored by at least three colors. Remark that this definition is considered only in the case where $d\geq 2$.
\end{defi}
Let us recall that by taking all faces of dimension $0$ and $1$ of a simplicial complex  $\bigtriangleup,$ i.e. all vertices and edges, we associate to $\bigtriangleup$ a simple graph $\Gcal_{\btu}$.
\begin{prop}\cite[Theorem 1.1]{bm1} Let $\btu$ be a simplicial complex of dimension $d$. Denote by $R_\btu$ the associated Stanley--Reisner ring. Assume that $\Gcal_{\btu}$ admits a good $(d+1)$--coloration. Let $\Ccal_1$, $\Ccal_2, \ldots, \Ccal_{d+1}$ sign the classes of colors and for each $ i = 1, 2, \ldots d+1,$ put
$$g_i = \sum_{x_j \in \Ccal_i}x_j.$$
Then $g_1,$ $g_2, \ldots, \ g_{d+1}$ is a system of parameters of $R_\btu.$ In particular, the reduction number of $R_\btu$ is 1.
\end{prop}

In \cite{bm2}, Barile and Morales considered a class of binomial ideals defining an union of scrolls, which indeed are binomial extension of  simplicial  ideals $\Bcal_{\ov \btu}$ where $\Gcal_{\btu}$ is a generalized $d$--tree and each vertex belongs to at most two extension facets. They proved that this class of ideals is of reduction number 1, and an explicit expression of the reduction is given. 

In this section, we will extend these results to the binomial extension of  simplicial  ideals.

\begin{note} Let $F_l = \left\{ x^{(l)}_0,x^{(l)}_1, \ldots, x^{(l)}_{i_1},\ldots,\ ,x^{(l)}_{i_{k_l}},...,x^{(l)}_{d_l} \right\}$ be an extension facet of $\btu$ by the points $\{y^{(l)}_{i_jm}\}$in the proper edges $\left\{ \ed{x^{(l)}_0}{x^{(l)}_{i_1}},\ldots,\ \ed{x^{(l)}_0}{x^{(l)}_{i_{k_l}}} \right\}$, and we denote $I_l=\left\{ x^{(l)}_0,x^{(l)}_{i_1},\ldots,\ ,x^{(l)}_{i_{k_l}} \right\}$. 
\end{note}

For the binomial extension of  simplicial  ideals, it is not necessary to color all $\Gcal_{\ov \btu}$. In fact, for each extension facet $F_l$, it is sufficient to color the extremal points in each bloc of the associated matrix $M_l$. The graph obtained from these points is defined as follows:

\begin{defi} The \emph {reduced graph}, denoted by $\wt \Gcal_{({\ov \btu}, \Bcal)}$, is given by:
\begin{itemize} \item The vertex set $V_{\wt \Gcal_{({\ov \btu},\Bcal)}}$ consists of the points of $\btu$ and the points $y^{(l)}_{i_j1}$ (with $j \ge 2$) for all extension facet $F_l$ of $\ov \btu$.
\item The set of edges is 
$$E_{\wt \Gcal_{({\ov \btu},\Bcal)}} := \bigcup_{F_l \in \ov \btu \rm{ s.t } Y^{(l)} \neq \emptyset} \left\{ {\left(E_\btu \setminus \{\ed{x^{(l)}_0}{x^{(l)}_{i_j}}\}^{k_l}_{j=2}\right)} \cup \{\ed{x^{(l)}_0}{y^{(l)}_{i_j1}}, \ed{x^{(l)}_{i_2}}{y^{(l)}_{i_j1}}\}^{k_l}_{j=2} \right\}.$$
\end{itemize}
\end{defi}

\begin{ex}\label{Greduit} Consider the following binomial extension with the binomial ideal associated to the facet $F=[a, b,c, d]$ extended by ${x, y, z}$ is generated by $2 \times 2$ minors of the matrix 
$$M = \left\{ \begin{array}{cc} a&x \\ x&b \end{array} \left| \begin{array}{c} y \\ c \end{array} \left| \begin{array}{c} z \\ d \end{array} \right. \right. \right\}.$$  
\begin{center}
\includegraphics[height=5.2cm,width=6.5cm]{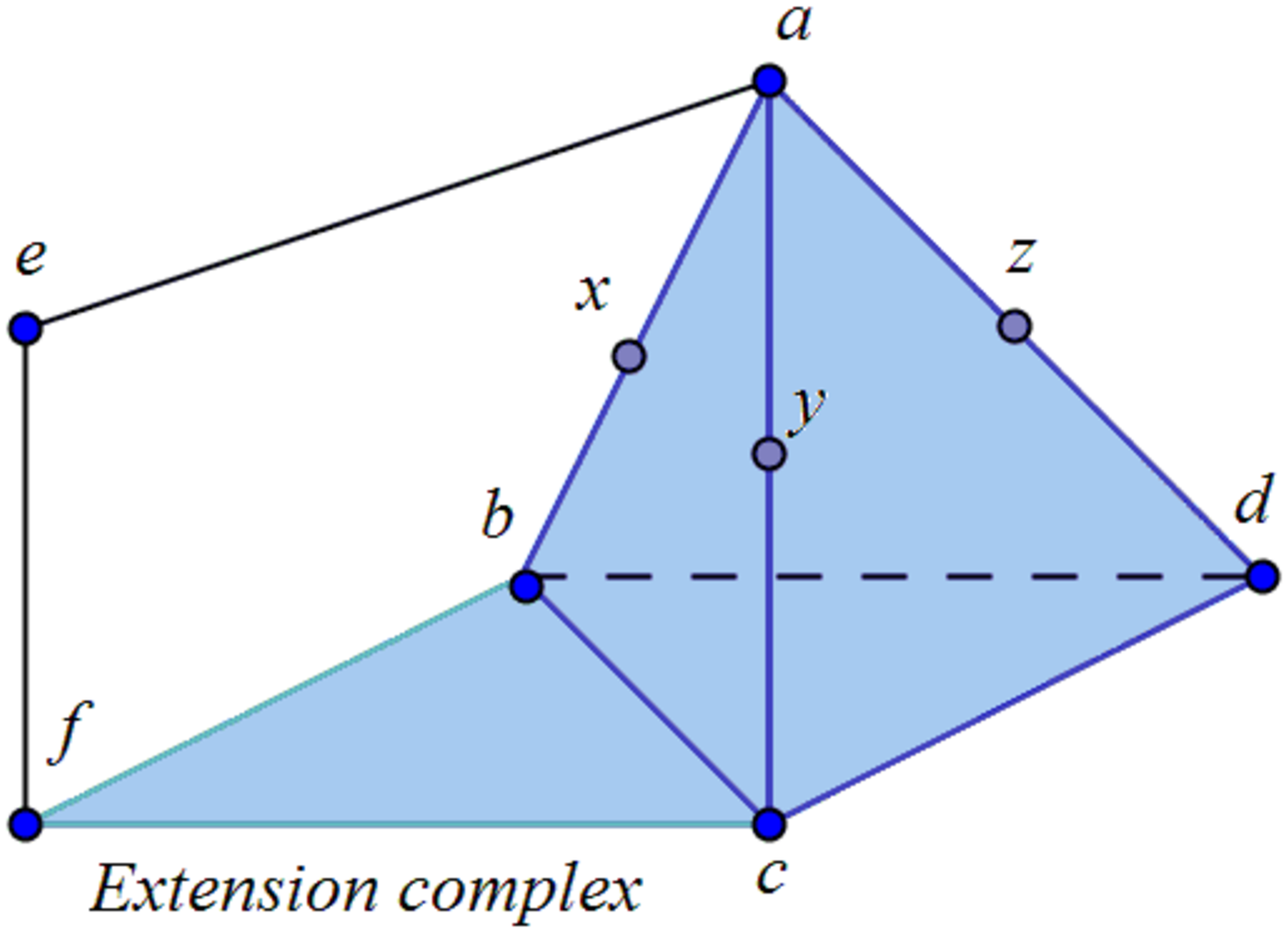} \includegraphics[height=5.2cm,width=6.5cm]{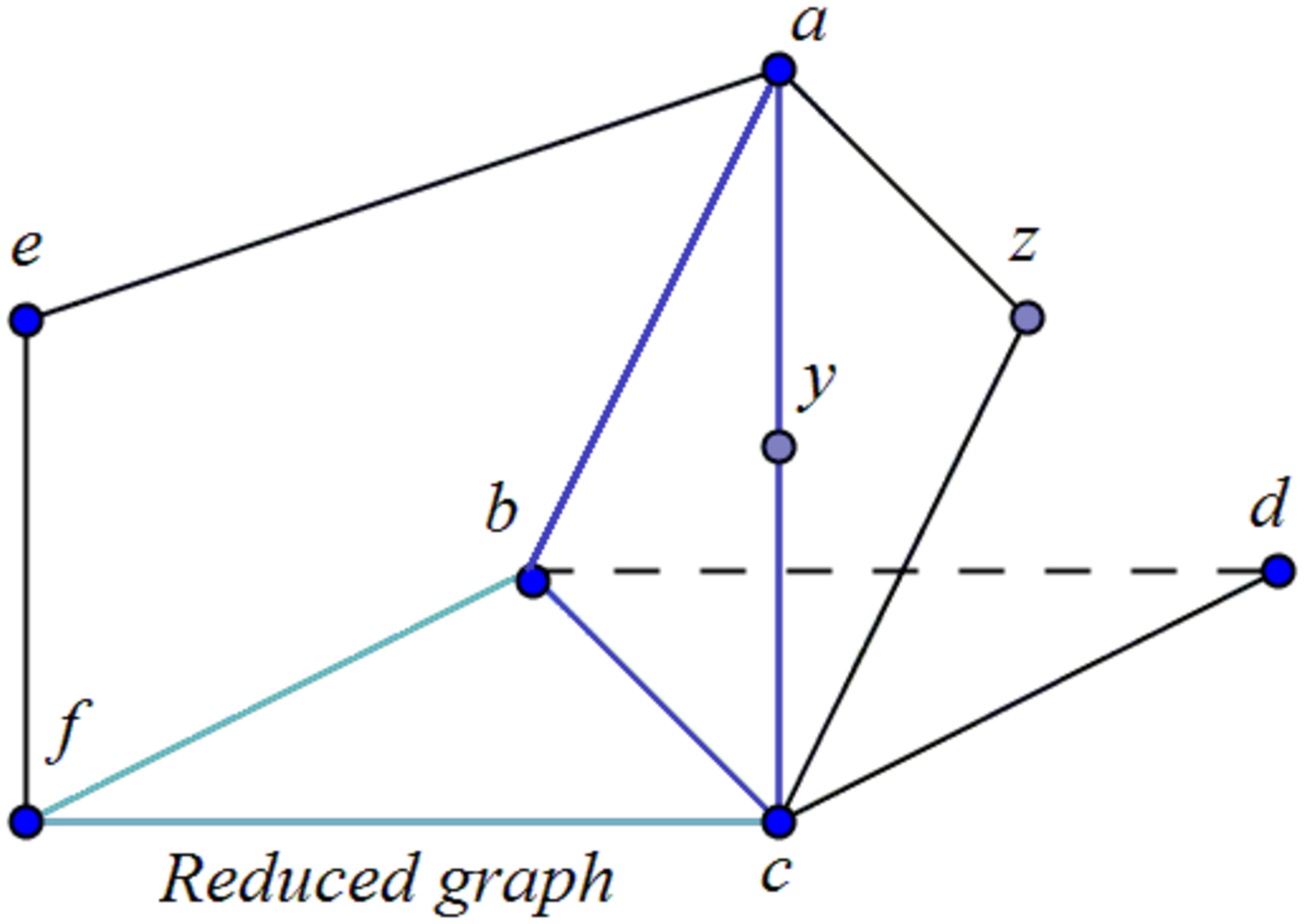}
\end{center}
Then, the reduced graph associated to this extension complex is as in the figure.
\end{ex}
\begin{lemma}\label{modB}
Let $\Bcal_M$ be the ideal generated by all the $2\times 2$ minors of the matrix $M$ 
$$ M := \left( \begin{array}{cccc} x_0 & y_{ 1,1} & \ldots & y_{ 1,j_{ 1}}\\ y_{ 1,1} 
& y_{ 1,2} & \ldots & x_{ 1} \end{array} \left|
\begin{array}{ccc} y_{ 2,1} & \ldots & y_{ 2,j_{ 2}} \\ y_{ 2,2} & \ldots & x_{ 2}
\end{array} \left|
\begin{array}{c} \ldots \\ \ldots \end{array} \left|
\begin{array}{ccc}  y_{ {k},1} & \ldots & y_{ {k},j_{ {k}}}\\ y_{ {k},2} &\ldots &
 x_{ {k}}
    \end{array}   \right. \right. \right.      \right).$$
and $\fr m'=\{x_0,...,x_k\}\cup Y$ be the set of variables in the matrix $M$,   then for any two distinct variables $x,y \in \fr m'$ with $y\in Y$ the product $xy$ is equivalent modulo $\Bcal_M$  to  one of the following monomials: 
\begin{enumerate}
\item $x_{ 1}x_0,$
\item $x_{ m}p, p\in Y_n$ 
\item $qy_{ n,1}, q\in Y_m,$ $n\geq 2$,
\item $x_0y_{ 1,j},$
\item $x_0y_{ n,1},$
\end{enumerate} 
Where $1\leq m\leq n\leq k$
\end{lemma}
Proof. We can assume that $x$ is a variable in the $m-$block of $M$, and $y$ is a variable in the $n-$block of $M$, with $m\leq n$. 
We have the following cases:
\begin{enumerate}
\item $1\leq m< n\leq k$
\begin{enumerate}
\item $x=x_{ m}$, \hfill\break In this case we have the monomials 
$x_{ m} p$,   with  $m\geq 1, p\in Y_n$,
\item $x$ any  and $y=y_{ n,1}$, \hfill\break In this case we have  the monomials 
$x_{0} y_{ n,1},q y_{ n,1}, x_{ m}y_{ n,1}$ with $q\in Y_m, m\geq 1$.
\item $x=y_{ m,u},y=y_{ n,v}$, or the case $x=y_{ m,u},y=x_{ n}$. We only consider the first case, the proof of the second case is similar.
$$xy=y_{ m,u}y_{ n,v}=(y_{ m,u}y_{ n,v}-y_{ m,u+1}y_{ n,v-1})+y_{ m,u+1}y_{ n,v-1}$$
$$=y_{ m,u+1}y_{ n,v-1} \mod \Bcal_M=...=\begin{cases}
x_{ m}p \mod \Bcal_M, p \in Y_n& \cr qy_{ n,1} \mod \Bcal_M, q \in Y_n& \end{cases}. $$
\end{enumerate}
\item $m=n>1$ 
\begin{enumerate}
\item $x=y_{ m,1}$.  In this case we have the monomials $y_{ m,1}p, y_{ m,1}x_{ m}$ with $p\in Y_m$.
\item $y=x_{ m}$. In this case we have the monomials $px_{ m}$ with $p\in Y_m$.
\item $x=y_{ m,u},y=y_{ m,v}$, with $u<v$. In this case we have 
$$xy=y_{ m,u}y_{ m,v}=(y_{ m,u}y_{ m,v}-y_{ m,u-1}y_{ m,v+1})+y_{ m,u-1}y_{ m,v+1}$$
$$=y_{ m,u-1}y_{ m,v+1} \mod \Bcal_M=...=\begin{cases}x_{ m}p \mod \Bcal_M, p \in Y_m& \cr qy_{ m,1} \mod \Bcal_M, q \in Y_m&\end{cases}. $$
\end{enumerate}
\item $m=n=1$
\begin{enumerate}
\item $x=x_{0}$. In this case we have the monomials $ x_{0}p$ with $p\in Y_1$.
\item $x=x_{1}$. In this case we have the monomials $ px_{ 1}$ with $p\in Y_1$.
\item $x=y_{ 1,u},y=y_{ 1,v}$. 
 $$xy=y_{ 1,u}y_{ 1,v}=(y_{ 1,u}y_{ 1,v}-y_{ 1,u-1}y_{ 1,v+1})+y_{ 1,u-1}y_{ 1,v+1}$$
$$=y_{ 1,u-1}y_{ 1,v+1} \mod \Bcal_M=...=\begin{cases}x_{ 1}p \mod \Bcal_M, p \in Y_1& \cr x_{ 1}p \mod \Bcal_M, p \in Y_1& \cr x_0x_{ 1} \mod \Bcal_M, p \in Y_1&\end{cases}. $$
\end{enumerate}
\end{enumerate}
\begin{defi} We will say that  $\wt \Gcal_{({\ov \btu}, \Bcal)}$ admits a binomial-coloration if  $\wt \Gcal_{({\ov \btu}, \Bcal)}$ admits a (d + 1)-coloration $\wt \Ccal$ such that  for every facet $F_l$:
\begin{enumerate}
\item $\wt \Ccal (x^{(l)}_0)\cap F_l = \{x^{(l)}_0,x^{(l)}_{ 2}\}$,
\item $\wt \Ccal (y^{(l)}_{ j1}) \cap F_l =\{y^{(l)}_{ j1},x^{(l)}_{ {j+1}}\} $ pour tout $(j = \ov{2, k_l-1}),$
\item  $\wt \Ccal(y_{ k,1}^{(l)})\cap F_l =\{y_{ k,1}^{(l)}\},$
\item For all the other vertices $x\in \wt \Gcal_{({\ov \btu}, \Bcal)}\cap F_l $, $\wt \Ccal (x)\cap F_l = \{x\} .$
\end{enumerate}
\end{defi}
 \begin{prop}\label{red0} 
 Suppose that $\wt \Gcal_{({\ov \btu}, \Bcal)}$ admits a binomial-coloration $\wt \Ccal$, we set 
$$g_i = \sum_{x \in \wt\Ccal_i}x,\ \ \ \Gbb:=(g_1,\ g_2, \ldots, \ g_{d+1}).$$
 Consider a facet  $F_l$ with a nonzero associated scroll matrix.
\begin{enumerate}
\item If $x^{(l)}_0y^{(l)}_{ n,1} \in \Gbb\fr m+\Bcal_{\overline\Delta}$ for any $2\leq n\leq k_l$, then $$xy\in \Gbb{\fr m}_{\ov \btu} + \Bcal_{\ov \btu} $$ for any variables $x\not=y\in F_l$, excepts for the products $x^{(l)}_0x^{(l)}_{1}, y^{(l)}_{ 1,u} y^{(l)}_{ 1,v}$, $x^{(l)}_0y^{(l)}$,  for $y^{(l)}\in F_l$ but not appearing in $M_l$, and $x^{(l)}y^{(l)}$,  for $x^{(l)},y^{(l)}\in F_l\cap \btu\setminus \{x^{(l)}_0\}.$  
\item If $x^{(l)}_0\in F_l^\circ ,$ then $$xy\in \Gbb{\fr m}_{\ov \btu} + \Bcal_{\ov \btu} $$ for any variables $x\not=y\in F_l$,  excepts for the products  $x^{(l)}y^{(l)}$,  for $x^{(l)},y^{(l)}\in F_l\cap \btu\setminus \{x^{(l)}_0\}.$   
\end{enumerate} 
\end{prop}{}

Proof.- We call this facet   $F$, we also delete all scripts $l$ from the variables defining vertices  in $F$ and the associated matrix $M$. from now on, we will denote by $\equiv $ the equivalence relation introduced by  $\Gbb\fr m+\Bcal_{\overline\Delta}$. We have the following: 
\begin{remark}\label{type1}\begin{enumerate}
\item If $x,y$ are two distinct element in $F$ such that $\Ccal(x)\cap F=\{x\},$ and $y$ belongs only to the facet $F$, 
(i.e. $y\in F^\circ $), then $xy\equiv  0$, since we can write
$$ xy=g_x y-\sum_{z\in C(x), z\not=x} zy,$$
but $y $ belongs only to the facet $F$ and $z\not\in F$, so $zy \equiv  0$, hence $xy\equiv  0$.
\item If $x,y$ are two distinct element in $F$ such that $\Ccal(x)\cap F=\{x,x'\},$ and $y\in F^\circ $ then 
$xy\equiv x'y$.
\end{enumerate}
\end{remark}
We have  the following cases:
\begin{enumerate}
\item $x$ doesn't appears in the matrix $M$ and $y\in Y$,
\item $x,y$ appear in the matrix $M$,  but one of them  belongs to $ Y$,
\item $x=x_0$, and   $y=x_m$, for some $1\leq m\leq  k$

\end{enumerate}

 Now we consider each case:
\begin{enumerate}
\item If $x$ doesn't appears in the matrix $M$ and $y\in Y$, then $\Ccal(x)\cap F=\{x\},$ the Remark \ref{type1} applies so  $xy\equiv  0$.\hfill\break
\item If $x,y$ appear in the matrix $M$,  but one of them  belongs to $ Y$.

By applying the Lemma \ref{modB}  the monomial $xy$ is equivalent modulo $ \Gbb m+\Bcal_{\overline\Delta}$ to one of the following monomials:
\begin{enumerate}
\item For  $1\leq m\leq n\leq k$, $x_{ m}p, p\in Y_n$,   
or $qy_{ n,1}, q\in Y_m,$ $n\geq 2$,
\begin{description}
\item[a-1) ] The remark \ref{type1} applies to the monomial  $x_{ 1}p, p\in Y_n, n\geq 1$, so it belongs to $\Gbb m+\Bcal_{\overline\Delta}$.
\item[a-2) ] If we are in the first case then by hypothesis $x_0y_{ n,1} \equiv  0$. If we are in the second case the following argument is also true when $q=x_0$. Consider the case $qy_{ n,1} , q\in Y_1 ,n\geq 2$.\hfill\break 
If  $n=k$,  and  $q\in Y_1$ the remark \ref{type1} applies, so $qy_{ k,1} \equiv  0$ for any $q\in Y\cup \{x_0\}, q\not=y_{ k,1}.$  So we may assume that $2\leq n<k$. In this case $C(y_{ n,1})\cap F= \{y_{ n,1}, x_{ n+1}\}$, by applying the remark \ref{type1} we have:
$$ qy_{ n,1}\equiv  -q x_{ {n+1}} ,$$ 
so by using the binomial relations in the matrix $M$ we have $$ qy_{ n,1}\equiv  -q x_{ {n+1}} \equiv \begin{cases}
x_{ 1}p' , p'\in Y_{n+1}& \cr
q'y_{ {n+1},1} , q'\in Y_1& \end{cases}.$$
Since the case $x_{ 1}p, p\in Y_n, n\geq 1$ was considered in the item b-1), after a finite number of steps we have either 
$qy_{ n,1} \equiv 0$ or $qy_{ n,1} \equiv  q' y_{ {k},1} $, but it was proved before that $q' y_{ {k},1}\equiv 0 $, so  $qy_{ n,1} \equiv 0 $.

\item[a-3) ] We now consider the monomial $x_{ 2}p$, with $p\in Y_n,$ $n\geq 2$. Since $C(x_{ 2})\cap F= \{x_{ 0}, x_{ 2}\}$, by applying the remark \ref{type1} we have  $$ x_{ 2}p\equiv -x_0 p.$$ 
By using the binomial relations in the matrix $M$ we have 
$$x_{0}p \equiv \begin{cases}x_{0}y_{ n,1} & \cr
x_{ 1}p' , p'\in Y_n& \cr
qy_{ n,1} , q\in Y_1& \end{cases},  $$
So this case is done by taking care of the previous cases. 
\item[a-4) ] We   consider the monomials $x_{ m}p$, for $2<m\leq n\leq k$, and the monomials $qy_{ {n},1}$ with $q\in Y_{m}$, $2\leq m<n\leq k$. Since $ \wt\Ccal(x_{ m}) \cap F)=\{x_{ m}, y_{ {m-1},1}\}$ by applying the remark \ref{type1} we have 
$$x_{ m}p\equiv  -y_{ {m-1},1}p .$$ By the proof of the Lemma \ref{modB}, case 1.c, we have 
$$y_{ {m-1},1}p \equiv \begin{cases}x_{ {m-1}}p & p\in Y_n  \cr qy_{ {n},1}
& q\in Y_{m-1}  \end{cases}.$$
We have either $$x_{ m}p\equiv x_{ {m-1}}p \equiv ...\equiv x_{ {2}}p , $$ or $$x_{ m}p\equiv qy_{ {n},1}
, $$
for some $q\in Y_{m}$, $2\leq m<n\leq k$. \hfill \break So it should be enough to consider the monomial 
  $qy_{ {n},1}$ with $q\in Y_{m}$, $2\leq m<n\leq k$. The case $n=k$  was considered in b-3). So we may assume $2\leq m<n< k$. By applying the remark \ref{type1} we have 
$$qy_{ {n},1}\equiv -qx_{ {n+1}} ,$$ and by using the binomial relations in $M$
$$qx_{ {n+1}} \equiv \begin{cases}x_{ m}p & p\in Y_{n+1}  \cr q'y_{ {n+1},1} &  q'\in Y_m\end{cases}.$$
So after a finite number of steps we will have: 
\begin{itemize}
\item Either $qy_{ {n},1}\equiv   x_{ 2}p $, $p\in Y_{n}, n\geq 3$, yet considered in b-1) or 
\item $qy_{ {n},1}\equiv   q'y_{ {k},1} $, $q'\in Y_m$, yet considered in b-2).
\end{itemize}
\end{description} 
\item If $x,y$ belongs to the first block of $M$ then we have to consider the following monomials $x_1y_{ 1,j},x_0y_{ 1,j},x_0x_1$ or $y_{ 1,u}y_{ 1,v}$. For the monomial $x_0y_{ 1,j}$, by applying the remark \ref{type1} we have  $$x_0y_{ 1,j}\equiv  -x_{2}y_{ 1,j} ,$$ which is equivalent modulo the binomial relations in the matrix $M$ either to the monomial $x_{ 1}p$, with $p\in Y_2$,  or to the monomial $qy_{ 2,1}$, with $q\in Y_1$. The monomial $x_{ 1}p$, with $p\in Y_2$, was yet considered in the first item, and the the monomial $qy_{ 2,1}$, with $q\in Y_1$ was considered in the second item.

The monomial $x_1y_{ 1,j}$, was consider  before. By the  Lemma \ref{modB} the monomial $y_{ 1,u}y_{ 1,v}$ is equivalent modulo $\Bcal_{\overline\Delta}$ to one of the monomials $x_1y_{ 1,j},x_0y_{ 1,j},x_0x_1$. If $x_0\in F^\circ $ then the remark \ref{type1} applies to $x_{ 1}x_0$ so we have $x_{ 1}x_0 \equiv  0$. So this subcase is done.
\end{enumerate}
\item We consider the monomial $x_0x_m$, where $x_m$ appears in the matrix $M$, $m>1$.
By the proof of the Lemma \ref{modB}, 1.c, for any $1<m$,   we have either $x_0x_{ {m}}\equiv x_{ 1}p$, with $p\in Y_{m}$, or 
$x_0x_{ {m}}\equiv qy_{ {m},1}$, with $q\in Y_{1} .$ Both monomials were considered in the previous items and we have seen that they belong to $ \Gbb m+\Bcal_{\overline\Delta}.$ 
The proposition is proved. \hfill $\Box$  
\end{enumerate}
\begin{prop}\label{Red1} 
 Suppose that $\wt \Gcal_{({\ov \btu}, \Bcal)}$ admits a binomial-coloration $\wt \Ccal$. Let  $\wt \Ccal_1$, $\wt \Ccal_2, \ldots, \wt \Ccal_{d+1}$ the classes of colors of the vertices of $\wt \Gcal_{({\ov \btu}, \Bcal)}$, we set 
$$g_i = \sum_{x \in \wt\Ccal_i}x,\ \ \ \Gbb:=(g_1,\ g_2, \ldots, \ g_{d+1}).$$
\begin{enumerate}
\item Suppose that  this  coloration on  the graph $\wt \Gcal'$, with set of edges $E_{\wt \Gcal'}:=E_{\wt \Gcal_{({\ov \btu},\Bcal)}}\cap E_{\Gcal_\Delta } $, is a good $d+1-$coloration $\wt \Ccal$, that is every cycle has more than three colors.  Let ${\fr m}_{\ov \btu} = \left({\rm {\bf x}}, {\rm {\bf y}} \right)$ be the maximal ideal of the  polynomial ring $\Rcal:= k[{\rm {\bf x}}, {\rm {\bf y}} ].$  
\item Suppose that   for every facet $F_l$ we have either:
\begin{enumerate}
\item $x^{(l)}_0y^{(l)}_{ n,1} \in \Gbb\fr m+\Bcal_{\overline\Delta}$ for any $2\leq n\leq k_l$.
\item $x^{(l)}_0\in F^\circ $.
\end{enumerate}
\end{enumerate}
Then we have :
$$\fr m^2_{\ov \btu}= \Gbb{\fr m}_{\ov \btu} + \Bcal_{\ov \btu} .$$
In particular, the reduction number of $\Rcal \slash \Bcal_{\ov \btu}$ is 1.
\end{prop}
 Proof.-   We consider a monomial $xy\in \fr m^2$. We want to prove that $xy\in \Gbb\fr m+\Bcal_{\overline\Delta}. $ 

$\bullet $ We have two  cases:

\begin{description}
\item[A.]  The variables $x,y$ are distinct. We remark that if $x,y\in\fr m$ are not in the same facet then $xy\in \Bcal_{\overline\Delta}$, hence  $xy\equiv 0$. \hfill\break
 $\bullet $  We can assume that the  variables $x,y$  belong to the same facet $F_l$. We fix this facet and we call it  $F$, we also delete all scripts $l$ from the variables defining vertices  in $F$ and the associated matrix $M$. 
By applying the Proposition \ref{red0}, we have to study the cases :$x=x_0$, and   $y=x_m$, for some $k< m\leq  d$ or $x=x_l$, and   $y=x_m$, for some $1\leq l< m\leq d$.
\item[B.]  The variables $x,y$ coincide.
\end{description}
We consider now in detail both cases.

{\bf A.}  
 If   $x_m$ doesn't appears in the matrix $M$ or $m=1$, then the edge $\la x_0,x_m\ra$ belongs to $\Gcal':=E_{\wt \Gcal_{({\ov \btu}, \Bcal)}} \cap E_{\btu}$. So it is enough to consider 	a monomial  $xy=x_l x_m$, such that $\la x_l,x_m\ra$    belongs to $\Gcal'$. In this case $C(x)\cap C(y)=\emptyset $ and the edge 
$\la x,y\ra$ belongs to  $\Gcal'$. We have  
$$xy= g_{x}y -\sum_{z\in C(x),z\not= x} zy ,$$ that is 
$$xy\equiv  -\sum_{z\in C(x), zy\not\in \Bcal_{\overline \Delta }} zy .$$
Let remark that the condition $zy\not\in \Bcal_{\overline \Delta }$ implies that $y, z$ belongs to the same facet and are distinct since they have distinct colors. Taking care of all solved cases we can assume that $z=x_s^{l'}, y=x_t^{l'}$ belongs to some facet $F_{l'}$ and the edge $\la z,y\ra$ belongs to $\Gcal'$.  By applying again the same argument we have:
$$xy\equiv  \sum_{z\in C(x), zy\not\in \Bcal_{\overline \Delta }} \sum_{w\in C(y), zw\not\in \Bcal_{\overline \Delta }}zw .$$
If it rests some monomials in the sum, we redo the algorithm. We will prove that the algorithm stops after  a finite number of steps. Assume the opposite that the algorithm will never stop. Then, there exists an infinite chain of edges $\la x,y\ra $, $\la y,z\ra $, $\la z,y _1\ra$, $\la y_1,z_1\ra $, \ldots. Since the number of the variables is finite, we must have a cycle in this chain, i.e. we have a cycle in  $\Gcal'$. Moreover, each point of this cycle is colored by either the color of $x$ or the one of $y$ (here, we have $\Ccal(x) \neq \Ccal(y)$ because $x$ and $y$ form an edge of  $\Gcal'$. This is a contradiction to the fact that $  \Gcal'$ admits a good coloration.

Hence, the algorithm will stop, i.e.:
$$xy \in (g_1,\ g_2, \ldots, \ g_{d+1}){\fr m}_{\ov \btu} + \Bcal_{\ov \btu}.$$
{\bf B.} To finish the proof, we consider the case $x = y.$ 

\noindent If $x$ is a vertex colored, then by replacing $x$ by $g_{x}$, one has
$$x^2 = xg_{x} - \sum_{z \in \Ccal(x),z\not=x} xz .$$
We have that  $z$ and $x$ have the same color, but they are distinct, so $xz\equiv  0$ by the case A). Hence, $x^2\equiv  0$.
\hfill\break If $x$ is not colored, then $x$ appears in some scroll matrix and we have a binomial  $x^2-yz\in \Bcal_{\overline \Delta }$, with $y\not= z$, so $x^2\equiv yz\equiv 0$, by the above cases. 
The proposition is proved. \hfill $\Box$  

\begin{ex}\label{Greduit1} Our Proposition applies to the example \ref{Greduit} colored as follows. 
\begin{center}
\includegraphics[height=5.2cm,width=6.5cm]{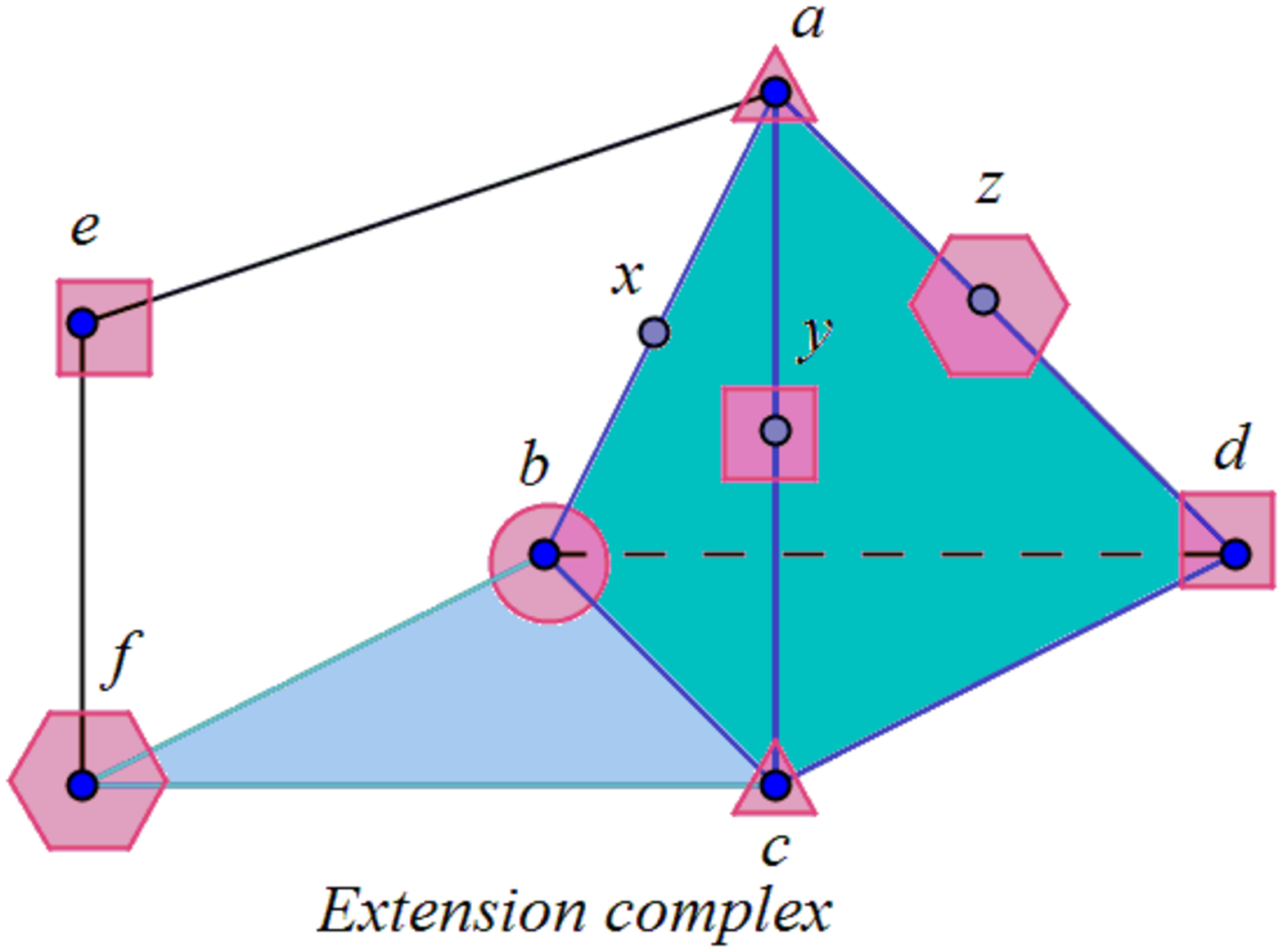} \includegraphics[height=5.2cm,width=6.5cm]{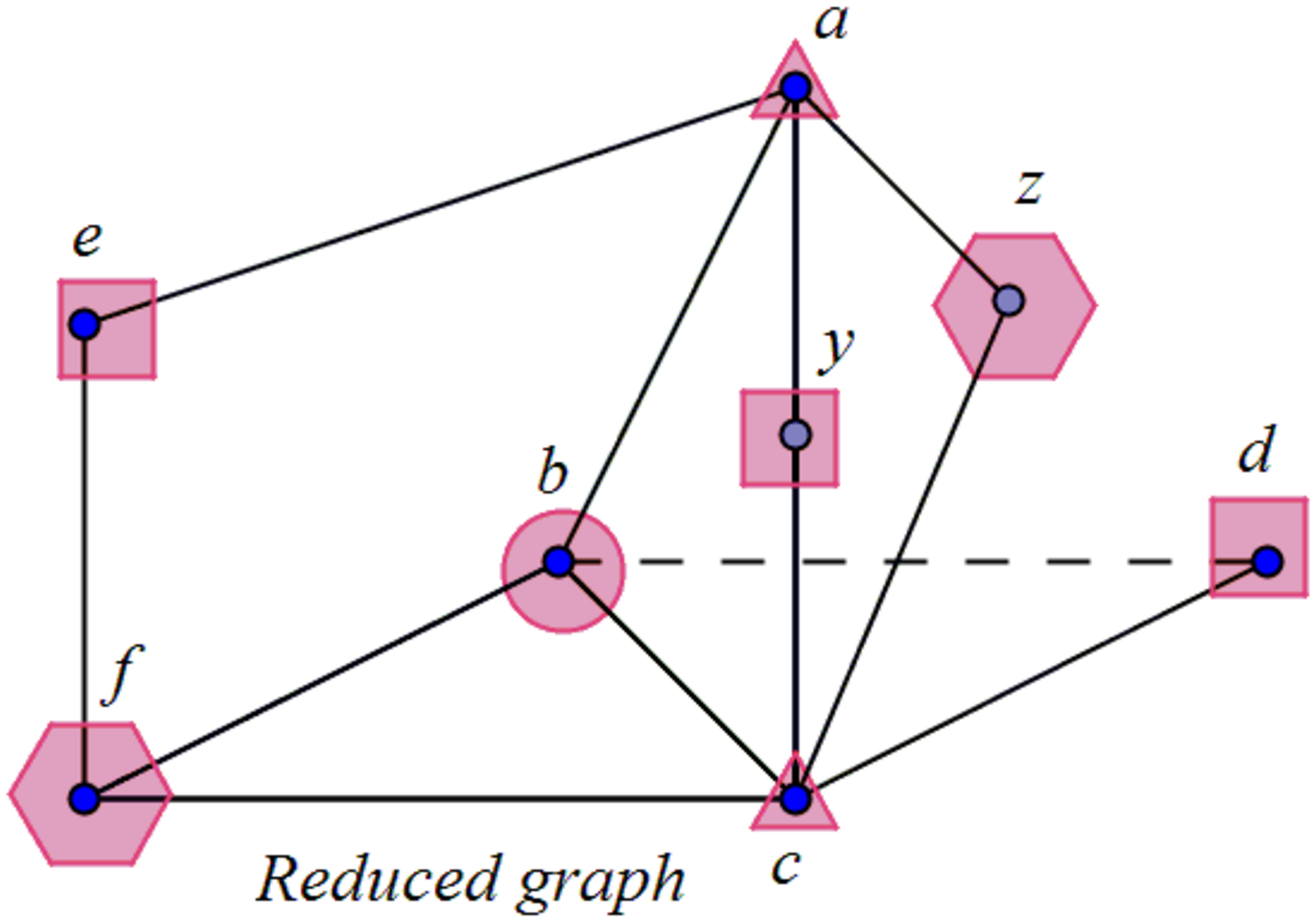}
\end{center}
\end{ex}
\begin{ex}\label{Greduit1} Our Proposition applies to the following complex \ref{Greduit1} colored as showed, and extended by the scroll matrices:
$$M_1:= \left( \begin{array}{cc} a\\b \end{array} \left| \begin{array}{c} y\\c \end{array} \right. \right), M_2:= \left( \begin{array}{cc} f\\c \end{array} \left| \begin{array}{c} z\\b \end{array} \right. \right), M_3:= \left( \begin{array}{cc} d\\f \end{array} \left| \begin{array}{c} x\\e \end{array} \right. \right), M_4:= \left( \begin{array}{cc} g\\a \end{array} \left| \begin{array}{c} w\\e \end{array} \right. \right)$$ 
\begin{center}
\includegraphics[height=5.2cm,width=6.5cm]{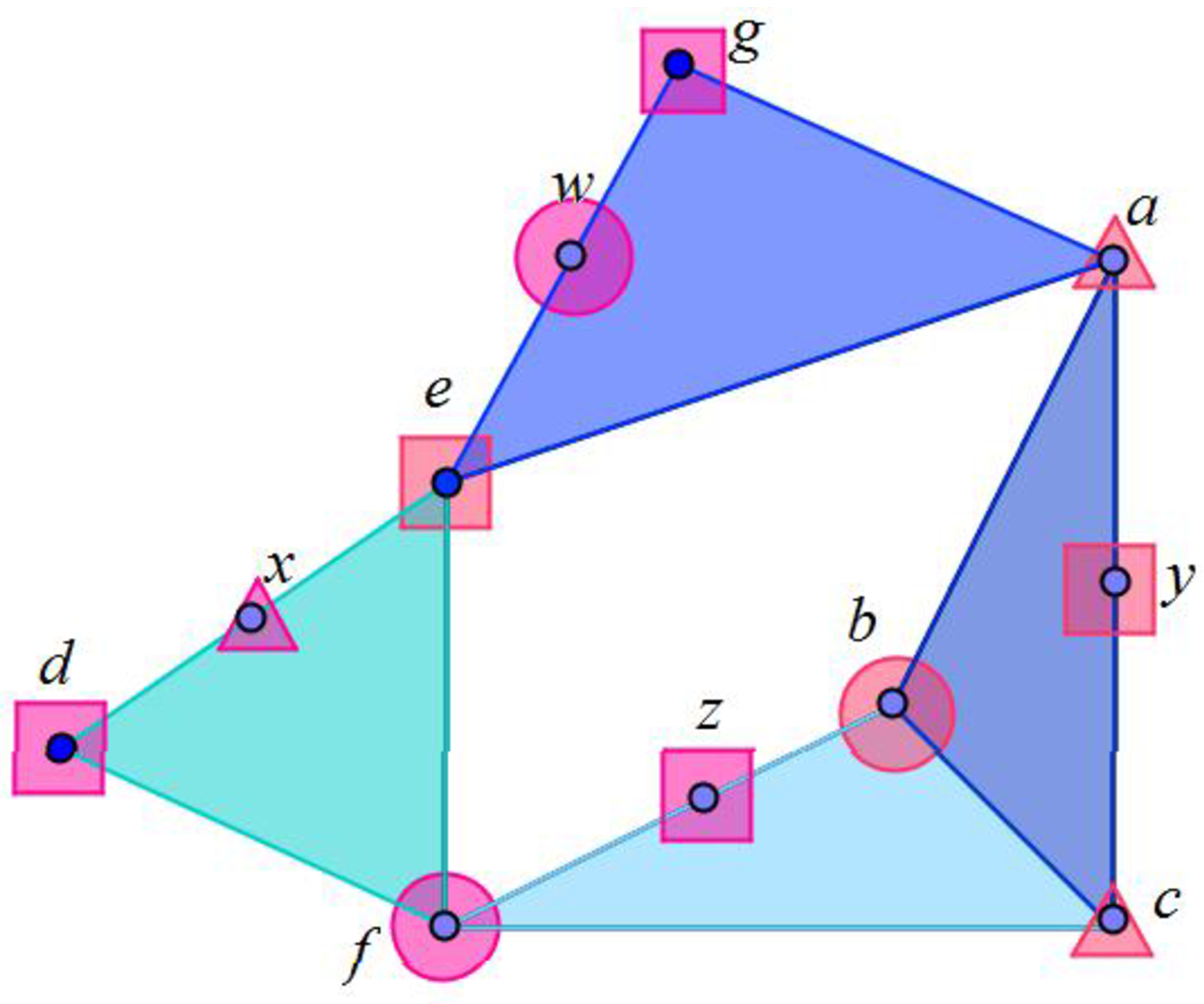} 
\end{center}
\end{ex}
\begin{ex}\label{cycles} Our Proposition  applies to  the complex represented by the  picture in the left, 
 but our Proposition cannot be applied to  the complex represented by the  picture in the right, in this case we have that the degree of the projective variety defined by this extended complex is 8, the codimension is 7. On the other hand by   a computation  we can check that the ideal $(a+c+d, b+e, f+v+y)$ is a reduction with  reduction number  two. The complex is extended by  the following matrices:
$$M_1:= \left( \begin{array}{cc} a&x\\x&b \end{array} \left| \begin{array}{c} y\\c \end{array} \right. \right)$$
associated to the facet $[\ a,\ b,\ c]$ extended by $x,\ y$ ;
$$M_2:= \left( \begin{array}{cc} d&u\\u&b \end{array} \left| \begin{array}{c} v\\c \end{array} \right. \right)$$
associated to the facet $[\ b,\ c,\ d]$ extended by $u,\ v$.
 \begin{center}
\includegraphics[height=5.2cm,width=5.5cm]{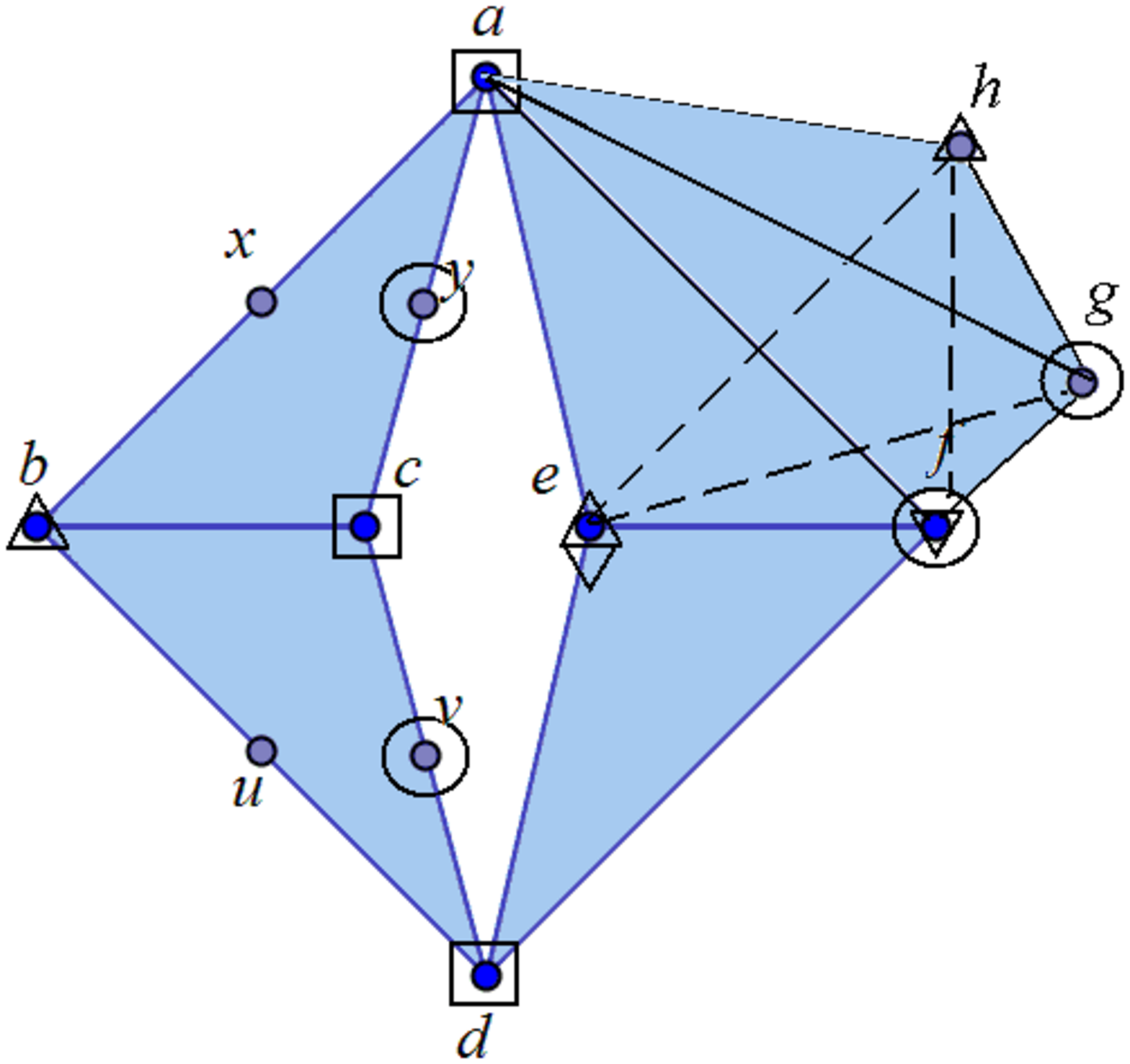}\includegraphics[height=5.2cm,width=5.5cm]{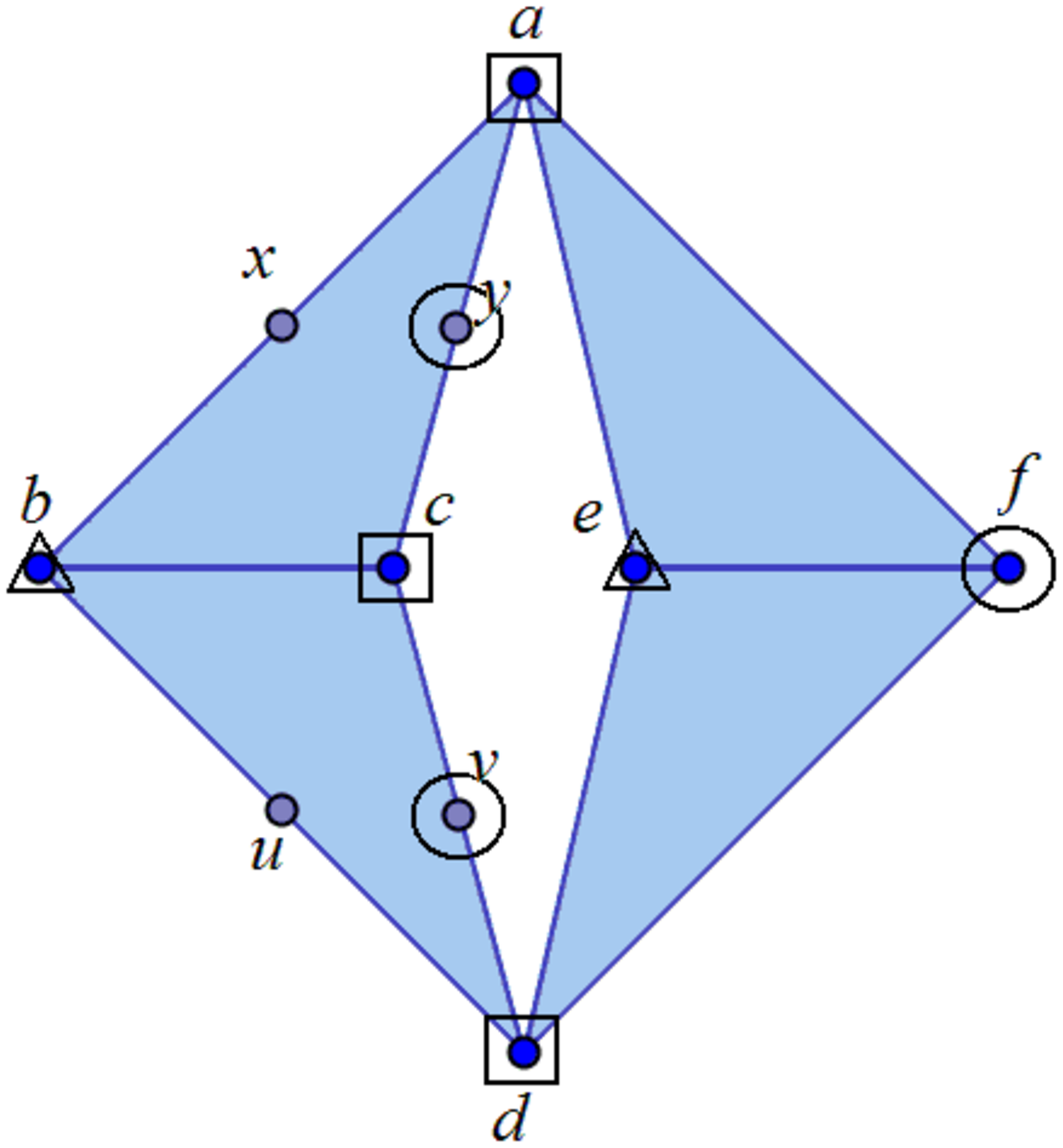}
\end{center} 
\end{ex}
 
If $\Gcal_{\btu}$ is a generalized 
$d$--tree, then in \cite{bm1} it was constructed an explicit reduction for the quotient by the Stanley-Reisner ideal associated to $\btu$.   The aim of the following proposition is to prove that the reduction number of $\Rcal \slash \Bcal_{\ov \btu}$ is 1, and to give an explicit expression of the reduction. 

\begin{prop}\label{Red2} Let $\Gcal_{\btu}$ be a generalized $d$--tree. Then $\wt \Gcal_{({\ov \btu},\Bcal)}$ 
admits a binomial--coloration. 
Let $\Ccal_1$, $\Ccal_2, \ldots, \Ccal_{d+1}$ denote the classes of colors. Put
$$g_i = \sum_{x \in \Ccal_i}x.$$
Then, we have:
$$(g_1,\ g_2, \ldots, \ g_{d+1}){\fr m}_{\ov \btu} + \Bcal_{\ov \btu} = \fr m^2_{\ov \btu},$$
where ${\fr m}_{\ov \btu} = \left({\rm {\bf x}}, {\rm {\bf y}} \right)$ is the irrelevant maximal ideal of the polynomial ring $\Rcal:= \Kcal[{\rm {\bf x}}, {\rm {\bf y}} ].$
\end{prop}

{\sc Proof:} The proof will be by  induction on the number $\lambda $ of facets of $\ov \btu$. 

\noindent
The case $\  \lambda =1$ is a particular case of the  Proposition \ref{Red1}. See also \cite{bm2}.

\noindent
Assume that the proposition is true for $\lambda  \geq 1$, we will prove it for $\lambda +1$. Since  $\Gcal_{\btu}$ 
is a generalized $d$--tree, one can find a facet $F$ such that its associated vertex in $H(\Gcal_{\btu})$ is a leaf. Consider $(\ov{\btu^\prime}, \Bcal^\prime)$ the extension complex constructed by the $\lambda $ facets different from $F$ in $\ov \btu$. We put $U = F \cap V_{\ov{\btu^\prime}}.$ The graph $\Gcal_{\btu^\prime}$ is also a generalized $d$--tree. By induction, the graph $\wt \Gcal_{(\ov {\btu^\prime},\Bcal^\prime)}$
 admits a good $(d+1)$--coloration as in the proposition, and for all $xy \in \fr m^2_{\ov {\btu^\prime}}$, we have:
$$xy \in (g^\prime_1,\ g^\prime_2, \ldots, \ g^\prime_{d+1}){\fr m}_{\ov {\btu^\prime}} + \Bcal_{\ov {\btu^\prime}},$$
where $\Ccal^\prime_i$, $g^\prime_i$ are the $i^{th}$ class of color and the correspondent sum. 
We have two cases:
\begin{description}
\item[I)] $F\cap  V_{\btu}=F$.
\item[II)] $F\cap  V_{\btu}\not=F$.
\end{description}
For each case, we will color the points in $F^\circ \cap V_{\wt \Gcal_{(\ov \btu, \Bcal)}}$, and we will define the sums $g_i.$ \\
\noindent
{\bf I):} $F\cap  V_{\btu}=F$. We can suppose that $F = \{x_0, x_1, \ x_2, \ldots, \ x_n \}$, and $U = \{x_1, \ x_2, \ldots, \ x_m \}$ ($1 \le m \le n \le d$), and $x_i \in \Ccal^\prime_i$ for all $i = \ov {1,m}.$ To obtain $(d+1)$--coloration which verifies the proposition, it is sufficient to color the points $x_i \notin U$ by arbitrary colors $\Ccal^\prime_j$ with $j > m.$
 For example, we can color $x_i$ by the $i^{th}$ color ($\forall i= \ov {m+1,n}),$ and $x_0$
 by the $(d+1)^{th}$ color. Hence, we have:
$$\Ccal_i = \Ccal^\prime_i \mbox{ and } g_i = g^\prime_i \mbox{ for all } i= \ov {1, m} \mbox{ or } i = \ov {n+1, d},$$
$$\Ccal_j = \Ccal^\prime_j \cup \{x_j\} \mbox{ and } g_j = g^\prime_j + x_j, j = \ov {m+1, n},$$
$$\Ccal_{d+1} = \Ccal^\prime_{d+1} \cup \{x_0\} \mbox{ and } g_{d+1} = g^\prime_{d+1} + x_0.$$

{\bf II):}  $F\cap V_{\btu}\not=F$. Let $M$ be the matrix associated to  $F$ :
$$M := \left( \begin{array}{cccc} x_0 & y_{1,1} & \ldots & y_{1,j_1}\\ y_{11} & y_{1,2} & \ldots & x_1 \end{array} \left|
\begin{array}{ccc} y_{2,1} & \ldots & y_{2,j_2} \\ y_{2,2} & \ldots & x_2
\end{array} \left|
\begin{array}{c} \ldots \\ \ldots \end{array} \left|
\begin{array}{ccc}  y_{l,1} & \ldots & y_{k,j_k}\\ y_{k,2} &\ldots & x_{k}    \end{array}   \right. \right. \right.      \right).$$
We can assume that $F \cap V_{\btu }= \{x_0, x_1, \ x_2, \ldots, \ x_l \}$ ($1 \le k \le l \le d$). Let remark that since by construction the propers edges  of $F$ are not in  $\ov{\btu^\prime}$, we have either $x_0 \in F^\circ$, or $x_i \in F^\circ$  $\forall i = \ov{1,k}.$ 
So we have to consider two sub-cases:

\begin{description}
\item[ II-1)] $x_0 \in F^\circ$, i.e. $x_0 \notin U$ : In order to color $\wt \Gcal_{({\ov \btu}, \Bcal)}$, we color each point $x\in  V_{\btu } \setminus \{U, x_0\}$ by a color not used in  $U$, and we define :
\begin{itemize} 
\item the color of  $x_0$ is the same of $x_2$ ;
\item the color of $y_{(j-1)1}$ is the color of $x_j$ for all $j = \ov{3,k}.$
\item the color of $y_{k,1}$ is the  $(l+1)^{\rm {i-th}}$ color not used in $F$. 
\end{itemize}
We can renumbering the classes of colors in such a way that $x_i \in \Ccal^\prime_i$ for $i = \ov{1,l}.$
Then, we have
$$\begin{array}{ccl} 
g_2& =&\begin{cases}
g^\prime_2 + x_0, \mbox{ if } x_{2} \in U,& \\ g^\prime_{2} + x_{2} + x_0, \mbox{ if } x_{2} \notin U; &\end{cases}\\ 
g_{j}& = &\begin{cases}g^\prime_{j} + y_{{(j-1)}1} \mbox{ for } j =\ov{3,l} \mbox{ such that } x_{j} \in U ,&\\ 
g^\prime_{j} + x_{j} + y_{{(j-1)}1} \mbox{ for } j =\ov{3,l} \mbox{ such that } x_{j} \notin U ,& \\ 
g^\prime_{j} + x_{j}\mbox{ for } j\in \{1,l+1,...,k\}, \mbox{ such that } x_{n} \notin U ;&\end{cases}\\ 
g_{l+1}&=&g^\prime_{l+1} + y_{k1};\\
g_j &=&g^\prime_j \mbox{ for all other indices $j$ ;}
\end{array} $$
% il faut colorier yl,1
\item[II-2)] If $x_0 \notin F^\circ$: In this case $x_i \in F^\circ$ for all $i = \ov{1,k}.$ We can suppose that $U = \{x_0, x_s, \ldots, x_l\}$ with $k < s \le l,$ and that $x_j \in \Ccal_j$ for all $j = \ov{s, l},$ and $x_0 \in \Ccal_2$. We put:
$$\begin{array}{lcl}
\Ccal_1 &= &\Ccal^\prime_1 \cup \{x_1\} ;\\
\Ccal_2 &= &\Ccal^\prime_2 \cup \{x_2\} ;\\
\Ccal_i &= &\Ccal^\prime_i \cup \{x_i, y_{(i-1)1}\} \mbox{ for all } i =\ov{3,k} ;\\
\Ccal_t &= &\Ccal^\prime_t \cup \{x_t\} \mbox{ for all } t =\ov{k+1,s-1};\\
\Ccal_j &= &\Ccal^\prime_j \mbox{ for all } j \ge s ; \\
\Ccal_{l+1} &= &\Ccal^\prime_{l+1} \cup \{y_{k,1}\}.
\end{array}$$
 Then we have:
$$\begin{array}{lcl}g_1 &= &g^\prime_1 + x_1;\\
g_2 &= &g^\prime_2 + x_2 ;\\
g_i &= &g^\prime_i + x_i+ y_{(i-1)1} \mbox{ for all } i =\ov{3,k} ;\\
g_t &= &g^\prime_t +x_t \mbox{ for all } t =\ov{l+1,s-1};\\
g_j &= &g^\prime_j \mbox{ for all } j \ge s \mbox{ and } j \neq k+1;\\
g_{l+1} &= &g^\prime_{l+1} + y_{k1}. \end{array}$$ 
\end{description}
Let us remark that in all the cases, for all $j$ the support of $g_j -g^\prime_j$ is contained in  $F^\circ$.

\begin{description}
\item[A) ] First we will prove that $xy\equiv 0$ for any $x\in F^\circ, x\not=y$ and $y\in F$.

$\bullet$ Case I)  We have that $\Ccal(y) \cap F$ contains only $y$, so by applying the Remark \ref{type1}, we have $xy\equiv 0$.

$\bullet$ Case II-1)  since $x_0 \in F^\circ, $ then by the Proposition \ref{red0} $xy\equiv 0$ for any variables $x\not=y\in F$,  excepts for the products  $xy$,  for $x,y\in F\cap \btu\setminus \{x_0\}$. Since we are interested in the monomials $x\in F^\circ, $ and $y\in F\cap \btu$, we have to consider the following cases:
\begin{description}
\item{II-1-a) }  $x \in F^\circ$, and $x= x_u,y=x_v$ appear in $M$, $u,v\not=0$. 
By applying the Remark \ref{type1}, since $x_u \in F^\circ$, we  have either $$x_u x_v\equiv 
\begin{cases}0 & {\rm if \ } v=1, \cr x_u x_0\equiv 0 & {\rm if \ } v=2, \cr x_uy_{(v-1)1}\equiv 0 & {\rm if \ } v >2. \cr
 \end{cases} $$
 \item{II-1-b)} $x \in F^\circ$. If  $x$ appears in $M$, and $y$ doesn't appears in $M$, since $\Ccal(y) \cap F = \{y\},$ we get by the   Remark \ref{type1},that $xy\equiv 0$. 
\item{II-1-c)} $x \in F^\circ$.  If $x$ doesn't appears in $M$ and $y$ appears in $M$, since $\Ccal(x) \cap F = \{x\},$ then due to the Remark \ref{type1}, it is sufficient to check the case where $y = x_v$ with $2 \le v \le k.$ But in this case, one has also that either
$  xx_v \equiv  xy_{(v-1)1}\equiv 0 $ or $xx_v \equiv  xx_0 \equiv 0$. 
\item{II-1-d)} $x \in F^\circ$, both $x,y$ don't appear in $M$:   One has $\Ccal(y) \cap F = \{y\},$ so $xy\equiv 0$ by to the Remark \ref{type1}. 
\end{description}
$\bullet$ Case II-2) First we prove that $x_0y_{n,1}\equiv 0$, for any $n\geq 2$. By the Remark \ref{type1}, we have 
$x_0y_{n,1}\equiv x_2y_{n,1}$. If $n=k$ since $x_2\in F^\circ , $ and  $\wt \Ccal (y_{k,1})\cap F = \{y_{k,1}\}$, Remark \ref{type1} we have $x_2y_{k,1}\equiv 0$. If $n<k$  then $x_2y_{n,1}\equiv  x_2x_{n+1}$, but $x_{n+1}\in F^\circ  $ so 
$x_2x_{n+1}\equiv  x_0x_{n+1}$. By using the binomial relations in the matrix $M$, we will have that 
$$x_0x_{n+1}\equiv \begin{cases}x_1 q& q\in Y_{n+1}\cr p y_{n+1,1}& p\in Y_1.\cr
\end{cases}$$
Now by the Remark \ref{type1} $x_1 q\equiv 0$, $ p y_{n+1,1} \equiv 0$ if $n+1=0$,  and $ p y_{n+1,1} \equiv p x_{n+1,1}$ if $n+1<k$. by applying the binomial relations in the matrix and the Remark \ref{type1}  after a finite number of steps we will have $x_0y_{n,1}\equiv 0$.

By using the Proposition \ref{red0}, we have  $xy\equiv 0$ for any variables $x\not=y\in F_l$, excepts for the products $x_0x_{1}, y_{ 1,u} y_{ 1,v}$, $x_0y$,  for $y\in F$ but not appearing in $M$, and  $xy$  for $x,y\in F_l\cap \btu\setminus \{x_0\}.$ So we need to consider the following cases: 
\begin{enumerate}
\item $x_0x_1\equiv -x_2x_1\equiv-x_2g{x_1}\equiv 0$. This case also will imply that 
\item $x_0 x_m$, $x_m\in F^\circ $ doesn't appears in $M$. $x_0x_m\equiv -x_2x_m\equiv-x_2g{x_m}\equiv 0$.
\item $x_m x_n, m,n>0$, $x_m\in F^\circ $, both  $x_m, x_n$ appear in $M$, this implies $x_n\in F^\circ $. We can assume that $n\geq 2$, so $x_mx_n\equiv x_m y_{n-1,1}$ if $n>2$, or $x_mx_n\equiv x_m x_0$, both monomial are equivalent to 0.
\item $x_m x_n, m,n>0$, $x_m\in F^\circ $, $x_m$ appears in $M$ but $x_n $ doesn't appears in $M$. Then  $x_mx_n\equiv-x_mg{x_n}\equiv 0$. 
\item $x_m x_n, m,n>0$, $x_m\in F^\circ $, $x_m$ doesn't appears in $M$ but $x_n $ appears in $M$. Then  $x_mx_n\equiv-g{x_m}x_n\equiv 0$. 
\item $x_m x_n, m,n>0$, $x_m\in F^\circ $, both $x_m,x_n$ don't appear in $M$. Then  $x_mx_n\equiv-x_mg{x_n}\equiv 0$. 
\end{enumerate}

\item[B) ]
$\bullet \ x , y\in V_{\ov{\btu^\prime}}$, By induction, one has
$xy = \sum_{i= 1}^{d+1} m_ig^\prime_i \ \mod \Bcal_{\ov{\btu^\prime}}$
with $m_i \in \fr m_{\ov{\btu^\prime}}.$
But
$$\sum_{i= 1}^{d+1} m_ig^\prime_i= \sum_{i= 1}^{d+1} m_i g_i -\sum_{i= 1}^{d+1} m_i(g_i-g^\prime_i).$$
Since the support of $g_i-g^\prime_i$ is in $F^\circ$ and Supp$(g_i-g^\prime_i)\cap$Supp$(m_i) = \emptyset$, due to the precedent cases 
$$m_i(g_i-g^\prime_i) \in (g_1,\ g_2, \ldots, \ g_{d+1}){\fr m}_{\ov {\btu}} + \Bcal_{\ov {\btu}}.$$ It implies that $xy$ verifies $(*).$  
 
\item[C) ]$\bullet \ x = y\in F^\circ$, In this case, if in addition $x$ is not colored, modulo $\Bcal_{\ov \btu}$ (see Lemma \ref{modB}), we re--obtain one of cases above. If $x$ is colored, we replace $x$ by the sum $g_x$ of all variables in the class of color of $x$, we will be in the case $x \neq y.$ Hence the proposition is proved. \hfill $\Box$  
\end{description}

\end{document}